\documentclass[12pt]{article}
\begin{document}
\noindent
\begin{center}
{\bf {\large A novel alternative to numerical integration }}\\    
\vskip 3mm
\noindent
{N. Mohankumar $^{a}$, Soubhadra Sen $^{b}$, A. Natarajan $^{c}$}   \\
\end{center}
\vskip 3mm
\noindent
{\it a.} 40, D J Nagar, Coimbatore, India 641004, Email: kovainmk@gmail.com\\
{\it b.} Radiological Safety Division, IGCAR, Kalpakkam, India 603102, Email: ssen@igcar.gov.in\\
{\it c.} 26, BBC City Park, II Phase, Chennai, India 600116, Email: anat1946@yahoo.co.in\\
\\
\begin{center}
{\bf Abstract}\\
\end{center}
\noindent
The Fundamental Theorem of Integral Calculus links the integrand and its antiderivative via a simple first order differential equation. A numerical solution of this ode yields the antiderivative and hence the required integral. This approach offers an economical and accurate alternative to the conventional approaches like the Gauss and the Double Exponential (DE) quadratures as demonstrated by a variety of examples.

\noindent
{\bf 1. Introduction}\\
\noindent
In physical sciences, quite often one encounters functions that can not be integrated in closed form and hence one has to resort to numerical integration. The Error function $Erf(x)$, the Exponential integral $Ei(x)$, and the Sine integral $Si(x)$ are some typical examples. For functions that can not be integrated exactly in closed form, one uses schemes like the Simpson quadrature or its generalisation, the Newton-Cotes method etc. If one needs better accuracy, usually the Gaussian quadrature is employed  [1,2]. If the function $f(x)$ is smooth then it is amenable to an accurate approximation by a polynomial and hence a Gauss quadrature of a particular order is guaranteed to give exact value. Thus for smooth functions, one normally prefers a Gaussian scheme.\\
\noindent
If the function has singularities like an end-point integrable singularity, a logarithmic singularity, or complex poles near the line of integration etc, then we set $\displaystyle f(x)\,=\,g(x)h(x)$. Here, $g(x)$ is the singular part of $f(x)$ and $h(x)$ is the smooth part. In this case, the quadrature sum is written as $\displaystyle \sum_{j=1}^n w_j\, h(x_j)  $ where the singular part $g(x)$ is absorbed in the weight function. This process of modifying the weights is an additional and rather involved task [3]. If the singularities depend on varying parameters, then the weights $\displaystyle\{w_j\}$ and the nodes $\displaystyle\{x_j\}$ do vary with these parameters which is a definite disadvantage. \\
\noindent
In subsequent discussions, an end-point integrable singularity, a logarithmic singularity or complex poles near the line of integration all will be referred to as singularities of {\it type A} for bevity. In recent times, a class of integration schemes like the TANH [4], the IMT [5] and the Double Exponential (DE) schemes [6,7] are widely employed. These schemes come under the category of {\it variable transformation method} and they are extremely handy for singularities of type {\it A}.  Essentially, all the three methods are just {\it simple} trapezoidal schemes {\it after} a specific change of integration variable. Here we outline the DE quadrature since it is superior to the other two methods. With $t$ and $u$ as the old and the new integration variables, the DE transformation introduced by Takahasi and Mori [6,7] is defined as follows.
\begin{eqnarray}
t &\in& (a,b); \, u \in (-\infty,\infty)\\
t_k \,&=&\,\phi (u_k)\,=\, (1/2)(b+a) + (1/2)(b-a)  \tanh(\frac {\pi}{2} \sinh(u_k))\\
\frac{dt}{du}\,&=&\,\phi'(u)\,=\,\frac{\pi (b-a)}{4} \, {\rm sech}^2  
[\frac{\pi}{2} \sinh(u)]  
\cosh(u)\\
u_k\,&=&\,kh, \, k=0,\pm1,\pm2,.....
\end{eqnarray}
$t_k$, the images of the equi-spaced nodes $u_k$ get clustered at the ends $t=a$ and $t=b$ of the old interval $(a,b)$. This specific property can be utilised to take care of the requirement of increased sampling near the poles of the integrand [8]. The derivative term $\displaystyle dt/du$ has the effect of annulling the integrable end-point singularity terms  like $1/(t-a)^\alpha$ or $1/(b-t)^\beta$ where $0 \,<\alpha,\beta\,<\,1$ and the logarithmic singularities like $\displaystyle log(t-a)$ or $\displaystyle log(b-t)$. This property is shared by the TANH, the IMT and the DE schemes [3] and hence they are well suited for handling singularities of type {\it A}.\\
\noindent
The Gaussian scheme is appropriate for smooth functions by its very nature of construction. But as we have seen earlier, when the integrand has singularities, {\it the method of modified moments} needs to be invoked which is  rather involved [3]. The TANH, the IMT and the DE methods are exceptionally well suited for handling integrals with type {\it A} singularities but they are not optimal for smooth integrands. In this context, it is worth exploring other possibilities. In the following we indicate {\it a powerful new approach based on approximating the antiderivative rather than the integrand} and it  is based on  the Fundamental Theorem of Integral Calculus.\\\\
\noindent
{\bf 2. Numerical Integration via the Fundamental Theorem of Integral Calculus}\\
\noindent
According to the Fundamental theorem of Integral Calculus [9], if $f(t)$ is continuous on the closed bounded interval $[a,b]$, and if $\displaystyle F(x)\,=\,\int_a^xf(t)\,dt,\,\,\,\, x\,\in[a,b] $, then $\displaystyle  \frac{dF}{dx}\,=\,f(x), \,\,x\in\,[a,b]$.
Given $f(t)$, this theorem indicates the possibility that if we solve the differential equation above for the antiderivative $F(x)$, then the required definite integral $\displaystyle \int_a^b f(t)\,dt  $   can be obtained  in terms of $F(x)$ as $\displaystyle  [F(b)\,-\,F(a)]$. As far as we know, not much has been explored to evaluate the required integral by this route even though this theorem is as old as calculus itself. Since we are dealing with a differential equation, we have a potential to {\it factor out the non-smooth or the dominant  part of the the antiderivative } which is a very big advantage. For smooth integrands, one can expand $F(x)$ in terms of a polynomial basis $\{\phi_n(x)\}$ as $\displaystyle  F(x)\,\simeq\,\sum _{j=1}^na_j\phi_j(x)$ and substitute this in the differential equation above and then find the expansion coefficients $\{a_j\}$ by a collocation process to finally evaluate $F(x)$. {\bf If the integrand has a dominant or a singular part, then the corresponding antiderivative component can be factored out whenever it is possible and again the solution can be attempted only for the smooth part of the antiderivative}. This results in a considerable reduction of the computational cost.\\
\noindent
Let the integrand $f(t)$ have  singularities like a pair of complex conjugate poles near the interval of integration. Then these poles will have a bearing on the  antiderivative function $F(x)$  and hence one may try a {\it rational approximation} in the following form.
\begin{eqnarray}
F(x)\,\simeq\, \frac {\sum _{j=0}^{n} a_jx^j}{1+\sum_{j=1}^m b_j x^j}
\end{eqnarray}
Then  fixing the $(n+m+1)$ unknowns involving the coefficients $\{a_j\}$ and $\{b_j\}$ by solving the differential equation by  collocation needs a non-linear iterative process. This appears to be computationally rather formidable and hence less desirable since the computational cost may turn out to be more expensive than that of our original numerical quadrature.\\\\
\noindent
{\bf 3. A Perturbation Approach for $F(x)$}\\
\noindent
In view of the difficulties mentioned above, we try  a perturbation approach for finding the  antiderivative F(x) that is best illustrated by a simple example.  Let $\displaystyle f(x)/x^{\beta},\, 0<\beta<1$ be the function that needs to be integrated over $(0,1]$. Here $f(x)$ is assumed to be smooth. The antiderivative of $\displaystyle 1/(x^{\beta})$ is $\displaystyle x^{1-\beta}/(1-\beta)$. Hence we try $\displaystyle  F(x)\,=\,\frac{\phi(x)\, x^{1-\beta}}{(1-\beta)}$. Here $\displaystyle \frac {x^{1-\beta}}{(1-\beta)}$ is the {\it " the dominant term "} and $\phi(x)$ is  {\it the "correction term"}  which is found out by solving the differential equation after substituting for $F(x)$ in terms of $\phi(x)$. Specifically we set $f(x)\,=\,e^{-x}$ and $\beta\,=\,1/2$. The required integral $\displaystyle \int_0^1
 \frac {e^{-x}}{\sqrt{x}}dx$, its alternate form (after a simple  transformation $x\,=\,t^2$)  and its exact value in terms of the Error function are given below.\\
\begin{eqnarray}
I\,=\,\int_0^1 \frac {e^{-x}}{\sqrt{x}}dx\,=\,2\int_0^1 e^{-t^2}dt\,=\,\sqrt{\pi}\,Erf(1)
\end{eqnarray}
Upon setting $\displaystyle F(x)\,=\, 2x^{1/2}\phi(x)$ the differential equation to be solved by collocation is given by
\begin{eqnarray}
2x \frac{d\phi}{dx}\,+\,\phi(x)\,=\,e^{-x};\,\, x\,\in[0,1]
\end{eqnarray}
In the first column of table 1, we give the results of our collocation solution as a function of the collocation nodes $n$. In columns two and three, we indicate the results of the Gauss-Legendre (GL) and the DE schemes as a function of the quadrature nodes $n$ for the unmodified integral. Since the function has an end-point integrable singularity, obviously the GL scheme converges poorly. The GL scheme with $250$ and $400$ nodes give the values $1.4901724491957$ and $1.4914742559640$, respectively and these values differ from the exact value $1.4936482656249$.  In all the tables $n$ indicates the quadrature nodes for the GL and the DE schemes while for the new method indicated here, $n$ stands for the collocation order. After the transformation $x=t^2$ the integral becomes $\displaystyle \int_0^1 2\,e^{-t^2}dt$ and now the integrand is smooth and has no singular component. Hence the GL scheme gives better results as indicated in table 2. Still, the new method lags only marginally behind the GL quadrature. However, for this smooth integrand, the DE scheme trails behind the other two methods and it needs  $47$ nodes to converge to $14$ digits.\\
\noindent 
The essence of our present approach can be stated succintly as,{\it " factor out the non-smooth or the dominant part of the antiderivative whenever possible and then fix the remaining part by collocation"}. In the subsequent sections, we give various examples illustrating the power of this approach and its superiority over both the Gaussian and DE methods.  Presently, we deal with quadrature problems in one dimension over a finite range with various singularities like poles, integrable singularities, removable discontinuities and non-integrable singularities lying just outside the range of integration. Other types of integrals will be dealt with in our subsequent papers. To solve the differential equation involving the correction function $\phi(x)$, we employ the Chebyshev pseudospectral method [10,11]. \\\\
\noindent
{\bf 4.1    Integrands with an end-point removable discontinuity}\\
\noindent
Next, we consider integrands which have a removable discontinuity at the ends of the interval like the following ones.\\
\begin{eqnarray}
I\,=\,\int_{0}^{1}\,\frac {sin(x)}{x}f(x)dx\, ;\,\,\, I\,=\,\int_{0}^{1}\,\frac {log(1+x)}{x}f(x)dx \nonumber
\end{eqnarray}
Due to the integrand term $(1/x)$, one may be tempted to write the antiderivative for the above integrals as $\displaystyle [log(x) \phi(x)]$ but this poses numerical problems. One can circumvent this difficulty in two ways. In the first method, the term $(1/x)$ is modified as  $\displaystyle (1/x^{\beta})$ where $0<\beta<1$ but $\beta$ is still very close to unity. Due to the restriction  $\beta<1$, the term $(1/x)$ is converted to an integrable weak singularity. Hence we seek the antiderivative as $\displaystyle [x^{1-\beta}\phi(x)/(1-\beta)]$ and then solve the resulting ode. As a consequence of setting $\beta$ very close but still not equal to unity, like $\beta\,=\, (1-10^{-16})$, a slight error is involved. However, this resulting error may be negligible. If more precision is needed, one can easily {\it extrapolate} the results for the case ($\displaystyle \beta \rightarrow \,1$), by a suitable sequence transformation like the Shanks transformation.\\
\noindent
In the second approach, instead of zero, we set the lower limit of integration as $\gamma$ where $\gamma \in [10^{-16},10^{-14}]$. Here $\gamma$ is very close but still not equal to zero.  As in the earlier approach, this modification results in an error which can be made negligible. We consider the test integral given below.
\begin{eqnarray}
I\,=\,\int_{0}^{1}\,\frac {log(1+x^2)}{x}dx
\end{eqnarray}
Following the first approach, we modify this integral as
\begin{eqnarray}
I\,=\,\int_{0}^{1}\,\frac {log(1+x^2)}{x^{\beta}}dx;\,\,\,\beta=(1-10^{-16})
\end{eqnarray}
By resorting to the second recipe, the required integral is modified as follows and it is solved by the pseudospectral method without any factorisation.
\begin{eqnarray}
I\,=\,\int_{\gamma}^{1}\,\frac {log(1+x^2)}{x}dx;\,\,\,\gamma=10^{-16}
\end{eqnarray}
Columns one and two of table $3$ provide the results of the new method employing the first and second approaches, respectively. Columns $3$ and $4$ are the results of the GL and the DE schemes, respectively for the unmodified integral. Once again, the GL method has an edge over the new method and the DE method has the least economy.\\\\
\noindent
{\bf 4.2  Singularity close to one end but just outside the integration interval}\\
\noindent
Next, we consider two  integrands where the singularity is {\it not} an integrable weak singularity but lies just outside the interval $[\gamma,1]$ at $x=0$. \\
{\bf Test problem 4.2.a}\\
The first problem in this category is given below.
\begin{eqnarray}
I\,=\,\int_{\gamma}^{1}\,\frac {e^{-x^2}}{x}dx\,;\,\,\,\,\, \gamma\,=\,10^{-5}
\end{eqnarray}
For the new and the GL methods, we change the integration variable by the rule $x\,=\,e^u$ and the transformed integral is given below.
\begin{eqnarray}
I\,=\,\int_{\log(\gamma)}^{0}\,e^{-(e^u)^2}du
\end{eqnarray}
After the transformation, the integrand is smooth. Hence for the new method, we resort to the pseudospectral solution of the ode relating to Eq.(12)  without any further factorisation of the antiderivative. For the DE method,  we evaluate both the integrals given in Eq.(11) and Eq.(12) and we confine to Eq.(12) for the GL evaluation. The results are indicated in table $4$. Here too the GL scheme scores over the other two methods. But the new method trails behind the GL scheme not by a big margin. It must be noted that if we calculate the integral of Eq.(11) by the GL method, then even with $400$ nodes only a two digit convergence is obtained.  Further, while employing the DE method for the integrals given in Eq. (11) and Eq. (12), we need $154$ and $74$ nodes, respectively, for $14$ digit convergence.\\
\noindent
{\bf Test problem 4.2.b}\\
We consider the integral given below.
\begin{eqnarray}
I\,=\,\int_{\gamma}^{1}\,\frac {log(1+e^{-x})}{x^{1.5}}dx\,;\,\,\,\,\,\, \gamma\,=\,10^{-5}
\end{eqnarray}
For the new method, we take the antiderivative as $\displaystyle [x^{-0.5}/(-0.5)]\,\phi(x)$ and then solve for $\phi(x)$ by collocation. The results are indicated in table $5$. In this case, the  DE scheme needs  $182$ nodes to give the converged value. Here, the GL method does not converge at all owing to the closeness of the non-integrable singularity even though it is not within the interval of integration. The superiority of the new method is obvious.\\\\
\noindent
{\bf  4.3 A non-factorable weak singularity that can be factored by a modification}\\
\noindent
Below, we consider an integrand that has an integrable singularity at $x=0$ that is close to one end $\gamma$ of the interval $[\gamma,1]$.
\begin{eqnarray}
I\,=\,\int_{\gamma}^{1}\,\frac {e^{-x^2}}{{sin^{\beta}(x)}}dx\,;\,\,\,\,\, \gamma=10^{-5};\,\,\,0<\beta<1
\end{eqnarray}
To facilitate factorisation, we multiply the numerator and the denominator by $x^{\beta}$ and this converts the integrand to the familiar form, $\displaystyle f(x)/x^{\beta}$ with $\displaystyle f(x)=\frac {{x}^{\beta}e^{-x^2}}{sin^{\beta}(x)}$.  Hence we have
\begin{eqnarray}
I\,=\,\int_{\gamma}^{1}\,\frac {e^{-x^2}}{{sin^{\beta}(x)}}dx\,=\,\int_{\gamma}^{1}\,\frac {f(x)}{{x}^{\beta}}dx\,;\,\,\,\,\, \gamma=10^{-5};\,\,\,\beta=0.7
\end{eqnarray}
The antiderivative for the integrand can be chosen as $\displaystyle [x^{1-\beta}\phi(x)]/(1-\beta)$ and we can solve for $\phi(x)$ by our collocation process and table $6$ provides the results. For convergence to $14$ digits, the DE scheme needs $144$ nodes while even with $n=200$, the GL method yields a value that is accurate to just two digits. The new method has a considerable edge over the other two methods.\\\\
\noindent
{\bf4.4 Integral with a logarithmic singularity} \\
\noindent
Let us consider the following problem with a  logarithmic singularity.
\begin{eqnarray}
I\,=\,\int_{0}^{1}\,log(x)\,f(x)\,dx
\end{eqnarray}
Here, $log(x)$ becomes unbounded as $\displaystyle x \rightarrow 0^+$. One may modify the integrand as $\displaystyle [x\, log(x) f(x)]/x$ and then factor out the term $x$ in the denominator. But this does not result in any significant computational gain. But a change of variable from $x$ to $u$ by the mapping  $x\,=\,e^{u}$ improves the situation. This is illustrated by the following integral with an exact value of $(-0.25)$.
\begin{eqnarray}
I\,=\,\int_{0}^{1}\,x\,log(x)\,dx
\end{eqnarray}
Table 7 shows the results.  For the GL and the DE methods, the integral is evaluated without {\it and} with the transformation $x\,=\,e^{u}$. The columns with the marking {\it trans} indicate the results after the change of variable from $x$ to $u$. For all the three methods, the $u$ interval is confined to $[log({10}^{-14}),1]$. With $n=200$, the GL method without transformation gives the value $(-0.2500000001554)$ and the DE method after the change of variable needs $n=87$ for convergence to exact value. The rapidity of convergence of the new method is quite impressive which is not matched by the other two methods.\\\\
\noindent
{\bf 4.5  An apparently non-factorable case}\\
\noindent
The following integral has a denominator term that appears like not factorable at first sight.
\begin{eqnarray}
I\,=\,\int_{0}^{1}\,\frac{log(1+x+e^{x})}{x^{0.2}+2}dx
\end{eqnarray}
However, we can try the following two transformations for our new scheme.\\\\
\noindent
Transformation {\bf A}:  We set $x^{0.2}=z$ and then factor out the term involving $z^{4}$ resulting from this substitution and this is followed by a change of variable from $z$ to $u$ given by $z+2\,=\,e^{u}$.\\\\
\noindent
Transformation {\bf B}:  We try the following change of variable.
\begin{eqnarray}
x^{0.2}+2\,&=&\,u^{2}+(\sqrt{2})^2\,;\,\,\,\,\,\,\,\,\,\ u\,=x^{(0.1)} \,\,\,\,\,\, 
\nonumber
\end{eqnarray}
After this,  due to the  denominator term $\displaystyle (u^2+2)$ of the resulting integrand, the antiderivative can be factored as $\displaystyle F(u)\,=\, tan^{-1}(u/\sqrt{2})\,\phi(u)$ and we solve for $\phi(u)$ by the pseudospectral method.\\
\noindent
The results are indicated in table $8$.  For the DE scheme, it is found that the original integral without any transformation converges faster and hence only those values are reported. On the other hand, with the substitutions $A$ and $B$, the convergence of the GL scheme is better. Here, we provide the GL results with the substitution $A$. It is to be noted that even with a $250$ order quadrature, the GL scheme without the above transformations converges to just $8$ digits. On the otherhand,  for full convergence the GL scheme needs just $n=16$ with the transformation $A$  and $n=18$ with the transformation $B$.  It is interesting to note that for the new scheme, the direct collocation solution of this integral without the above transformations will be accurate to just about $6$ digits even with $n$ as large as $200$. The GL method has a slight edge over the new method for this test integral.\\\\
\noindent
{\bf 5. Test problems from published literature}\\
\noindent
In this section, we indicate the results of few test problems from the published literature.\\\\
\noindent
{\bf 5.1 Problem from  Hasegawa}\\
The following integral is from Hasegawa [12].
\begin{eqnarray}
\nonumber
I\,=\,\int_{-1}^{1}\,\frac{1}{(x^2+1)(x-c)^2}dx\,=\,5.000000099582060E8;\,\,\,\,c\,=\,-1-(1.E-09)
\end{eqnarray}
In the new method, we factor out  the terms corresponding to the two factors of the denominator of the integrand above. The denominator term $(x-c)^2$ results in a non-integrable singularity at $x=c$ that is too close to $(-1)$, the lower end of the integration interval. The results are given in table $9$ and the computational gain of the new method is obvious. Here the GL method does not converge at all. With increasing $n$, both the new and the DE schemes do not give any improvement beyond the relative error of the order of $10^{-8}$.\\\\
\noindent
{\bf 5.2 Problem from Mori and Sugihara}\\ 
\noindent
The following integral $\displaystyle \int_0^1 log(1/x)\,dx/x^{0.25}$ \,from  Mori and Sugihara [13] has an exact value $1.7777777777778 $. The integrand has both a logarithmic as well as an integrable end-point singularity at the lower end $x=0$. In the new scheme, first we factor out the term corresponding to $(1/x^{0.25})$. This is followed by a change of variable from $x$ to $u$ by the mapping  $x\,=\,e^u$.  Here $u$ is restricted to the range $[log(10^{-16}),0]$ for all the three methods while for the integration without the transformation $x \in (0,1)$. For comparison,  we indicate the relative errors of the new method, the DE and the GL  schemes with and without this transformation in table $10$. Here the term  'trans' within the bracket indicates the results after the above mentioned transformation. With $n=58$, the DE scheme with the transformation has an error  $(2.824E-11)$ and any further increase in $n$ does not yield an improvement. Even with $n=200$, the GL scheme without the transformation still has a large error of value $(8.363E-04)$. The rate of convergence of the new method is impressive. The relative error of the new scheme reaches a saturation very quickly and it does not decrease further beyond $(2.864E-11)$. This is due to the restriction of the lower limit in the $u$ variable to $log(10^{-16})$. A way to cut this error further by extrapolation is discussed in the appendix. \\\\
\noindent
{\bf 5.3 Test problem from Lether involving complex poles}\\
\noindent
In the following we consider an integrand  like $\displaystyle f(x)/[(x-a)^2+b^2]$ that has a pair of complex conjugate poles located at $\displaystyle (a\pm ib)$, close to the integration interval. In our method, the denominator of this integrand namely
$\displaystyle [(x-a)^2+b^2]$ can be considered as the dominant term and the antiderivative can be taken as $\displaystyle tan^{-1} [(x-a)/b]\,\phi(x)$ and as usual, $\phi(x)$ is found by the pseudospectral method. Specifically we consider the following test case  $\displaystyle \int_{-1}^{1}\,{e^{x}}dx/{[(x)^2+{(0.01)}^2]}$ from Lether [14] with a reference value $313.17205623933$.  Here, one may factor out the term involving  ($e^x$) also. For this problem it is observed that this does not have much effect on the convergence. Hence we report the values without factoring this term. The relative errors are indicated in table $11$. With $n=200$, the relative errors for the GL and the DE schemes are found to be  $(1.703E-13)$ and $(8.294E-10)$, respectively. Here, the new method scores over the other two methods since we have the facility to factor out the troublesome pole related term.\\\\
\noindent
{\bf 5.4 Test problem from Bailey}\\
\noindent
The following very interesting test integral is from Bailey [15], $\displaystyle \int_{0}^{\pi/2}\,\sqrt{tan(x)}\,dx\,=\, (\pi/\sqrt{2})$. By simple manipulations we can modify the integral as follows.
\begin{eqnarray}
I\,=\,\int_{0}^{\pi/2}\,\sqrt{tan(x)}dx\,=\,\int_{0}^{\pi/4}\,\left[\frac{sin(x)+cos(x)}{\sqrt{x}}\right]\,\sqrt{\frac{2x}{sin(2x)}}dx
\end{eqnarray}
\indent
Now in the modified integral above, one can factor out the term corresponding to $\sqrt{x}$ in the denominator from the antiderivative and the remaining part can be evaluated by the pseudospectral method. Also, due to the numerical precision related problems associated with the term $\displaystyle [(2x)/sin(2x)]$ at $x=0$, for the new method alone the range is chosen as $[10^{-16},\pi/4]$ instead of $[0, \pi/4]$ and the results are provided in table $12$. The new scheme reaches its saturation accuracy with just $n=12$ and it beats the other two schemes convincingly. The error saturation is due to the restriction of the lower limit of integration to $10^{-16}$ instead of zero. This error can be cut down further by extrapolation techniques. With $n=100$, the GL method gives a relative error of the order of $10^{-3}$ for both the original and the modified integrals and the values reported under the GL heading are for the modified integral. The  DE scheme  for the modified integral needs more nodes than the new method but it converges ultimately to  full $14$ digits when $n=36$. The evaluation of the original integral by the DE scheme results in an accuracy saturating at $10^{-9}$ and then onwards oscillations set in.\\\\
\noindent
{\bf 6. Discussions}\\
\indent
If the function to be integrated is smooth, then the Gaussian scheme is obviously the method of choice. On the otherhand, if the function has singularities of type {\it A}, then one will prefer the DE scheme or its variants  over the Gaussian schemes. However, if the function under consideration does not belong to the above indicated types, it is worth seeking other alternatives. The present method based on approximating the antiderivative rather than the integrand is one such attempt. The power of the new method stems from the following two key factors. Firstly, whenever possible, we can factor out the dominant or the non-smooth part of the antiderivative which significantly reduces the cost of approximation. Secondly, the use of the pseudospectral method  enhances the accuracy and the economy.  The new method has a considerable edge over the GL and the DE methods in many test cases considered here. Further, even in the test cases where the GL or the DE methods dominate, the present method lags {\it only marginally} behind. In all the comparisons, if a variable transformation exists for  the GL and the DE methods, then only the relatively better results from the evaluation of the modified {\it or} the unmodified integrals is reported. This makes the comparison fair. Here, we have not considered a very exhaustive set of various types of test integrals like the semi-infinite or the full infinte range integrals, trigonometric integrals etc. Also, we have not addressed  the convergence related questions. These aspects will be explored in future works. Finally, the new method has ample scope for the user to invent ingenious ways of factoring out the dominant/singular part of the antiderivative which is quite interesting. Except for the theoretical aspects of the pseudospectral solution, the required manipulations like factoring out the dominant part are well within the reach of a beginner that includes the high school student who has learnt calculus!\\\\
\noindent
{\bf 7. Conclusions}\\
\noindent
We have indicated a new method of numerical integration that is based on approximating the antiderivative that essentially involves the pseudospectral solution of a first order ode. This new method has a considerable advantage in several cases and it has a potential both to supplement and surpass the existing classical methods.\\\\
{\bf Acknowledgement}\\
The authors wish to acknowledge the influence of  C. Schwartz. His novel TANH method gave the clue to look beyond the Gaussian quadrature schemes.\\\\
\noindent
{\bf Appendix}\\
\noindent
{\bf Shanks transformation for extrapolation}\\
In the following, we briefly indicate the Shanks transformation applied to a test integral result. If we have a triplet, $s(n), s(n+m)$ and $s(n+2m)$ which
are equispaced members of a sequence, then the {\it extrapolated limit} of this sequence $\displaystyle S_{\infty}$ can be evaluated as follows [16].
\begin{eqnarray}
a&=&s(n);\,\,b=s(n+m);\,\,c=s(n+2m) \nonumber\\
S_{\infty}\,&=&\, \frac {(ac-b^2)}{(a+c-2b)}\,=\,a-\,\frac{(b-a)^2}{(a+c-2b)} \nonumber
\end{eqnarray}
To minimise numerical errors, it is preferable to use the second relation above for $\displaystyle S_{\infty}$.
After applying this transformation to a given set of members of a sequence, we will have a new set of once transformed quantities. We can apply Shanks transform to every consecutive triplet of this new set and get another one where each member is a result of a double Shanks transformation. This process can be repeated till convergence is obtained.\\
\noindent
The integral $\displaystyle \int_0^1 log(1/x)\,dx/x^{0.25}$ \,from  Mori and Sugihara [13] has an exact value $1.7777777777778 $. 
For this problem the relative error of the new scheme reaches a saturation very quickly and it does not decrease further beyond $(2.864)\, 10^{-11}$. This is due to the restriction of the lower limit of integration to the quantity  $\gamma=log(10^{-16})$. A way to cut the error further by extrapolation is discussed here. For a sequence of  $\gamma$ values, $\displaystyle  log(10^{-10}),log(10^{-11}),\dots,log(10^{-16})$, we evaluate this integral  with $n=20$ by our new scheme. These integral values are subjected to a Shanks transformation and the transformation is repeated. It is seen from the last column of table $13$ that after two  successive Shanks transformations, the maximum relative error of the order of $10^{-7}$ is reduced to about $10^{-14}$.\\\\
\noindent
{\bf Table 1}\\
Results for the integral  $\displaystyle \int_0^1 \frac {e^{-x}}{\sqrt{x}}dx\,=\,1.4936482656249$\\
\begin{tabular}{c c c c}
\hline
n&New&GL&DE\\
\hline
5&1.4936485952414&1.3346090298085&1.5927631143854\\
6&1.4936482588994&1.3592563874781&1.4209378459472\\
7&1.4936482657617&1.3772716644230&1.5309835005364\\
8&1.4936482656222&1.3910193370197&1.4720325318878\\
9&1.4936482656249&1.4018574919485&1.5058148394444\\
10&1.4936482656249&1.4106225005522&1.4880045124688\\
35&1.4936482656249&1.4691191401786&1.4936482656240\\
38&1.4936482656249&1.4710308658017&1.4936482656249\\
\hline
\end{tabular}
\\\\\\\\
\noindent
{\bf Table 2}\\
 Results for the non-singular integral $\displaystyle \int_0^1 2\,e^{-t^2}dt\,=\,1.4936482656249$\\
\begin{tabular}{c c c c}
\hline
n&New&GL&DE\\
\hline
6&1.4946289062500&1.4936482657803&1.7927165840294\\
7&1.4936475753784&1.4936482656233&1.3518534632933\\
8&1.4936518669128&1.4936482656249&1.5629883064697\\
9&1.4936484098434&1.4936482656249&1.4592324712248\\
15&1.4936482656244&1.4936482656249&1.4931996510320\\
16&1.4936482656249&1.4936482656249&1.4938537847693\\
\hline
\end{tabular}
\newpage
\noindent
{\bf Table 3 }\\
Results for the integral $\displaystyle\int_{0}^{1}\,\frac {log(1+x^2)}{x}dx$ with a converged value $0.4112335167121$\\
\begin{tabular}{c c c c c}
\hline
n&New (method I)&New (method II)&GL&DE\\
\hline
4&0.4100875586003&0.4096679687500&0.4112340814293&0.7448862919169\\
9&0.4112336011915&0.4112336635590&0.4112335167121&0.3957330733966\\
12&0.4112335166910&0.4112335166865&0.4112335167121&0.4137489421612\\
16&0.4112335167103&0.4112335167092&0.4112335167121&0.4114231645971\\
18&0.4112335167121&0.4112335167120&0.4112335167121&0.4112827175014\\
19&0.4112335167121&0.4112335167121&0.4112335167121&0.4112088017243\\
25&0.4112335167121&0.4112335167121&0.4112335167121&0.4112332215339\\
45&0.4112335167121&0.4112335167121&0.4112335167121&0.4112335167132\\
49&0.4112335167121&0.4112335167121&0.4112335167121&0.4112335167121\\
\hline
\end{tabular}
\\\\\\\\
\noindent
{\bf Table 4}\\ Results for the integrals $\displaystyle \int_{\gamma}^{1}\,\frac {e^{-x^2}}{x}dx\,=\,\displaystyle\int_{\log(\gamma)}^{0}\,e^{-(e^u)^2}du$ with $\gamma\,=\,\,10^{-5} $ and having a converged value $11.114625665372$.\\\\
\begin{tabular}{c c c c c}
\hline
n&New&GL&DE(Eq. 11)&DE(Eq. 12)\\
\hline
8&11.128906250000&11.114336218795&9.3306760027737&11.270589778037\\
24&11.114624977112&11.114625665370&10.985845884430&11.114090373227\\
26&11.114625665359&11.114625665372&11.202921694065&11.115039213196\\
40&11.114625665134&11.114625665372&11.110676417713&11.114628215472\\
45&11.114625665372&11.114625665372&11.116799502110&11.114625944461\\
50&11.114625665372&11.114625665372&11.114709630403&11.114625682226\\
60&11.114625665372&11.114625665372&11.114692422004&11.114625664886\\
\hline
\end{tabular}
\newpage
\noindent
{\bf Table 5}\\ Results for the integral $\displaystyle\int_{\gamma}^{1}\,\frac {log(1+e^{-x})}{x^{1.5}}dx$ with $ \gamma\,=\,10^{-5} $ and having a converged value $436.08354183864$.\\\\
\begin{tabular}{c c c}
\hline
n&New&DE\\
\hline
4&436.08312877543&1063.5154981277\\
8&436.08354183491&558.88403384888\\
9&436.08354183809&121.24780208925\\
10&436.08354183865&327.07094429924\\
12&436.08354183864&430.75863332422\\
14&436.08354183864&431.74536268890\\
40&436.08354183864&433.80030057011\\
60&436.08354183864&436.07265948137\\
\hline
\end{tabular}
\\\\\\\\
\noindent
{\bf Table 6}\\
 Results for the integral $\displaystyle \int_{\gamma}^{1}\,\frac {e^{-x^2}}{{sin^{\beta}(x)}}dx$ with $ \gamma=10^{-5}$ and $\beta=0.7$ and  having a converged value $2.9195592066824$.\\\\
\begin{tabular}{c c c c}
\hline
n&New&GL&DE\\
\hline
5&2.9196220802647&2.2401359196449&2.9668943734058\\
13&2.9195592066825&2.5664170875376&2.9183729939852\\
14&2.9195592066824&2.5855625264712&2.9404764251821\\
15&2.9195592066824&2.6026960287340&2.9212114974862\\
20&2.9195592066824&2.6673955198716&2.9131595567403\\
80&2.9195592066824&2.8602330274062&2.9195592444079\\
90&2.9195592066824&2.8695426408396&2.9195592043411\\
\hline
\end{tabular}
\newpage
\noindent
{\bf Table 7}\\ 
Results for the integral $\displaystyle\int_{0}^{1}\,x\,log(x)\,dx $ with an exact value $(-0.25)$.\\
\begin{tabular}{c c c c c c}
\hline
n&New&GL&GL(trans)&DE&DE(trans)\\
\hline
2&-0.25&-0.2578533682318&-0.0001328959735&-1.1568426104738&-0.0000000000173\\
3&-0.25&-0.2518457084908&-0.0227336228624&-0.0006945844762&-0.0879872452593\\
10&-0.25&-0.2500208786906&-0.2631142225538&-0.2630380257838&-0.3885709040680\\
14&-0.25&-0.2500057123735&-0.2501176579217&-0.2505449891516&-0.1525194883310\\
18&-0.25&-0.2500021512010&-0.2500001723709&-0.2500153026394&-0.3029588884950\\
22&-0.25&-0.2500009821174&-0.2500000000576&-0.2500003410778&-0.2311761479918\\
26&-0.25&-0.2500005100748&-0.2500000000000&-0.2500000065405&-0.2536181385900\\
30&-0.25&-0.2500002905650&-0.2500000000000&-0.2500000001130&-0.2502839472564\\
34&-0.25&-0.2500001774388&-0.2500000000000&-0.2500000000018&-0.2495707915608\\
38&-0.25&-0.2500001143930&-0.2500000000000&-0.2500000000000&-0.2501252768546\\
\hline
\end{tabular}
\\\\\\\\
\noindent
{\bf Table 8}\\
 Results for the integral $\displaystyle \int_{0}^{1}\,\frac{log(1+x+e^{x})}{x^{0.2}+2}dx$ with a converged value $0.3988251952415$.\\
\begin{tabular}{c c c c c}
\hline
n&New(A)&New(B)&GL&DE\\
\hline
12&0.3988252029349&0.3988235407263&0.3988251952360&0.3991827611351\\
14&0.3988251953449&0.3988250207832&0.3988251952414&0.3988804348125\\
16&0.3988251952561&0.3988251830555&0.3988251952415&0.3988332640645\\
18&0.3988251952441&0.3988251954659&0.3988251952415&0.3988260668126\\
20&0.3988251952414&0.3988251953721&0.3988251952415&0.3988252520998\\
23&0.3988251952415&0.3988251952445&0.3988251952415&0.3988252016271\\
27&0.3988251952415&0.3988251952415&0.3988251952415&0.3988251956235\\
35&0.3988251952415&0.3988251952415&0.3988251952415&0.3988251952415\\
\hline
\end{tabular}

\newpage
\noindent
{\bf Table 9}\\
Relative errors for the integral$\displaystyle\int_{-1}^{1}\,\frac{1}{(x^2+1)(x-c)^2}dx $\\
\begin{tabular}{c c c}
\hline
n&New&DE\\
\hline
2&1.022E-07&1.000E-00\\
5&9.984E-08&4.389E-00\\
10&9.822E-08&1.694E-00\\
20&9.673E-08&3.496E-01\\
50&9.483E-08& 9.444E-02\\
110&9.323E-08& 3.783E-04\\
150&9.261E-08&4.763E-06\\
200&9.203E-08&6.769E-08\\
\hline
\end{tabular}
\\\\\\\\
\noindent
{\bf Table 10}\\
 Relative errors for the integral  $\displaystyle \int_0^1 log(1/x)\,dx/x^{0.25}$.\\
\begin{tabular}{c c c c c c}
\hline
n&New&GL&GL(trans)&DE&DE(trans)\\
\hline
2&2.864E-11&2.844E-01&7.651E-01&5.477E-01&9.987E-01\\
10&2.864E-11&4.495E-02&8.612E-05&4.043E-03&6.268E-01\\
14&2.864E-11&2.937E-02&3.573E-09&7.310e-05&3.297E-02\\
18&2.864E-11&2.126E-02&2.863E-11&1.210E-06&6.924E-02\\
35&2.865E-11&8.878E-03&2.864E-11&3.372E-15&5.319E-05\\
26&2.864E-11&1.315E-02&2.864E-11&2.006E-10&9.537E-04\\
40&2.864E-11&7.433E-03&2.864e-11&1.249E-14&3.936E-06\\
\hline
\end{tabular}

\newpage
\noindent
{\bf Table 11}\\
Relative errors for the integral $\displaystyle \int_{-1}^{1}\,{e^{x}}dx/{[(x)^2+{(0.01)}^2]}$\\
\begin{tabular}{c c c c}
\hline
n&New&GL&DE\\
\hline

10&2.430E-05&6.286E-01&8.492E-02\\
20&9.271E-06&5.701E-02&4.049E-01\\
40&4.533E-07&5.906E-03&3.540E-02\\
80& 3.035E-10&3.728E-05&3.951E-04\\
100&3.944E-11&2.398E-06&1.010E-04\\
120&5.683E-12&1.391E-07&1.624E-05\\
140&1.303E-13&7.253E-09&2.007E-06\\
\hline
\end{tabular}
\\\\\\\\
\noindent
{\bf Table 12}\\
 Relative errors for the integral $\displaystyle\int_{0}^{\pi/4}\,\left[\frac{sin(x)+cos(x)}{\sqrt{x}}\right]\,\sqrt{\frac{2x}{sin(2x)}}\,dx $ \\\\
\begin{tabular}{c c c c}
\hline
n&New&GL&DE\\
\hline

4&2.832E-05&7.689E-02&3.140E-01\\
8&1.024E-08&4.082E-02&1.206E-02\\
12&9.003E-09&2.778E-02&2.629E-04\\
20&9.003E-09&1.694E-02&1.419E-09\\
30&9.003E-09&1.139E-02&2.010E-11\\
36&9.003E-09&9.516E-03&3.339E-14\\
\hline
\end{tabular}\\
\newpage
\noindent
{\bf Table 13\\}
Relative errors with Shanks transformation for the integral $\displaystyle \int_0^1 log(1/x)\,dx/x^{0.25}$ by the new method\\
\begin{tabular}{c c c c}
\hline
Lower Limit $\gamma$&relative error&Shanks1&Shanks2\\
\hline
$\displaystyle log(10^{-10})$&5.777E-07&....&.....\\
$\displaystyle log(10^{-11})$&1.125E-07&...&.....\\
$\displaystyle log(10^{-12})$&2.172E-08&2.518E-10&.....\\
$\displaystyle log(10^{-13})$&4.170E-09&4.075E-11&.....\\
$\displaystyle log(10^{-14})$&7.962E-10&6.640E-12&6.510E-14\\
$\displaystyle log(10^{-15})$&1.513E-10&1.082E-12&7.485E-16\\
$\displaystyle log(10^{-16})$&2.864E-11&1.681E-13&1.148E-14\\
\hline
\end{tabular}
\\\\\\\\
\noindent
{\bf References:}\\
1. Davis PJ, Rabinowitz P. 1984 {\it Methods of Numerical Integration}. London, UK: Academic Press.\\
2. Froberg CE. 1969 {\it Introduction to Numerical Analysis}. London, UK: Addison-Wesley.\\ 
3. Press WH, Teukolsky SA, Vetterling WT, Flannery BP. 2007 {\it Numerical Recipes}. New York, NY: Cambridge University Press.\\
4. Schwartz C. 1969 Numerical Integration of Analytic Functions. {\it J. Comput. Phys}. {\bf 4}, 19-29. (doi:10.1016/0021-9991(69)90037-0)\\
5. Iri M, Moriguti S, Takasawa Y. 1987 On a Certain Quadrature Formula. {\it J. Comput. Appl. Math.} {\bf 17}, 3-20. (doi:10.1016/0377-0427(87)90034-3)\\
6. Takahasi H, Mori M. 1974 Double Exponential Formulas for Numerical Integration. {\it Publ. RIMS.} {\bf 9}, 721-741.\\
7. Mori M. 1991 Developments in the Double Exponential Formulas for Numerical Integration. In {\it Intern. Congr. Math. 1990 (ICM-90), koyto, Japan, 21-29 August}, pp. 1585-1594. Springer-Verlag.\\
8. Natarajan A, Mohankumar N. 1993 On the numerical evaluation of the generalised Fermi-Dirac integrals. {\it Comput. Phys. Commun.} {\bf 76}, 48-50. (doi:10.1016/0010-4655(93)90118-V)\\
9. Goldberg RR. 1964 {\ itMethods of Real Analysis}. Massachusetts, USA: Blaisdell.\\
10. Canuto C, Hussaini MY, Quarteroni A, Zang TA. 2006 {\it Spectral Methods}. Berlin, Genramy: Springer-Verlag.\\
11. Javidi M, Golbabai A. 2007 Spectral collocation method for parabolic partial differential equations with Neumann boundary conditions. {\it Appl. Math. Sci.} {\bf 1}, 211-218.\\
12. Hasegawa T. 1997 Numerical integration of functions with poles near the interval of integration. {\it J. Comput. Appl. Math.} {\bf 87}, 339-359. (doi:10.1016/S0377-0427(97)00197-0)\\
13. Mori M, Sugihara M. 2001 The double-exponential transformation in numerical analysis. {\it J. Comput. Appl. Math.} {\bf 127}, 287-296. (doi:10.1016/S0377-0427(00)00501-X)\\
14. Lether FG. 1997 Modified quadrature formulas for functions with nearby poles. {\it J. Comput. Appl. Math.} {\bf 3}, 3-9.(doi:10.1016/0771-050X(77)90017-1)\\
15. Bailey DH, Jeyabalan K, Li XS. 2005 A Comparison of Three High-Precision Quadrature Schemes. {\it Exper. Math.} {\bf 14}, 317-329. (doi:10.1080/10586458.2005.10128931)\\
16. Brezinski C, Zaglia MR. 1991 {\it Extrapolation Methods: Theory and Practice}. Amsterdam, Netherland: Elsevier.
\end{document}